\documentclass{amsart}

\usepackage{amssymb}

\theoremstyle{plain}

\numberwithin{equation}{section}

\newtheorem*{teorema}{Theorem A}
\newtheorem*{teoremab}{Theorem B}

\newtheorem{cora}[equation]{Corollary}

\newtheorem{hyp}[equation]{Hypotheses}

\newtheorem{lem}[equation]{Lemma}

\newtheorem{prop}[equation]{Proposition}

\newtheorem{claim}[equation]{Claim}

\newcommand{\Irr}{\operatorname{Irr}}
\newcommand{\Lin}{\operatorname{Lin}}

\newcommand{\modu}{\operatorname{mod}}

\newcommand{\Ker}{\operatorname{Ker}}

\theoremstyle{definition}

\newtheorem{definition}[equation]{Definition}

\begin{document}	
\title{Products of characters and finite $p$-groups II}

\author{Edith Adan-Bante}

\address{Department of Mathematics, 
University of Illinois at Urbana-Champaign,
Urbana, 61801}

\email{adanbant@math.uiuc.edu}


\subjclass{20c15}

\date{2003}

\begin{abstract} Let $p$ be a prime number.
 Let $G$ be a finite  $p$-group and $\chi\in \Irr(G)$.
 Denote by $\overline{\chi}\in \Irr(G)$ the complex
conjugate of $\chi$. Assume that 
$\chi(1)=p^n$.  We show that the number of distinct irreducible 
constituents of the product 
$\chi\overline{\chi}$ is at least $2n(p-1)+1$.
\end{abstract}

\maketitle

\begin{section}{Introduction}
Let $G$ be a finite group. 
Denote by $\Irr(G)$ the set of 
irreducible complex characters of $G$. 
Let $\chi \in \Irr(G)$. Define $\overline{\chi}(g)$ to be
the complex conjugate $\overline{\chi(g)}$ of $\chi(g)$ 
for all $g \in G$. Observe that $\overline{\chi}\in \Irr(G)$.
 Through this work, we will use
the notation of \cite{isaacs}.


Let $\chi, \psi \in \Irr(G)$.
 Since the product of 
characters is a character, $\chi\psi$ is a 
character of $G$. So it can be expressed as an integral linear
combination of irreducible characters. 
Let
 $\eta(\chi, \psi)$ be the number of distinct irreducible
constituents of the product  $\chi\psi$.

The main result is the following

\begin{teorema}
Let $p$ be a prime number. Let $G$ be a finite $p$-group and 
$\chi\in \Irr(G)$ with $\chi(1)=p^n$.  Then 
$$\eta(\chi, \overline{\chi})\geq 2n(p-1)+1.$$
\end{teorema}

In Proposition \ref{example} we show that, given a prime
$p$ and an integer $n\geq 0$,  there exist a $p$-group
$G$ and a character $\chi\in \Irr(G)$ such that
$\chi(1)=p^n$ and $\eta(\chi, \overline{\chi})=2n(p-1)+1$. Thus the 
bound in Theorem A is optimal. 

Fix an odd prime $p$. In \cite{me}, we considered 
the set of possible values
that $\eta(\chi, \psi)$ can have for any finite $p$-group
$G$ and faithful characters $\chi, \psi \in \Irr(G)$. 
There we proved that there is a ``gap'' among the
set of all possible values that $\eta(\chi, \psi)$ can have, namely
either $\eta(\chi,\psi)=1$ or $\eta(\chi, \psi)\geq \frac{p+1}{2}$.
The following theorem is another example of ``gaps''.

\begin{teoremab}
Let $p$ be a prime number. Let $G$ be a finite $p$-group and 
$\chi\in \Irr(G)$.
Then one of the following holds:

(i) $\chi(1)=1$ and $\eta(\chi,\overline{\chi})=1$.

(ii) $\chi(1)=p$ and $\eta(\chi,\overline{\chi})=2p-1$ or
$p^2$.

(iii)  $\chi(1)\geq p^2$ and  $\eta(\chi,\overline{\chi})\geq 4p-3$.

\end{teoremab}

We can check that $p^2\geq 4p-3$ if  $p>2$. Thus 
 Theorem B implies that
 if $\eta(\chi,\overline{\chi})< 4p-3$ and $p>2$, then either
 $\eta(\chi, \overline{\chi})=1$ or 
$\eta(\chi, \overline{\chi})= 2p-1$.
\end{section}

\begin{section}{Notation}

Given a  character   $\Theta$ of $G$,
we denote by $\Xi (\Theta)=\{\chi\in \Irr(G)\mid [\chi, \Theta]>0\}$ 
the set of irreducible constituents 
of $\Theta$. 

Given subgroups $N$ and $M$ of $G$, and characters
 $\nu\in \Irr(N)$ and $\mu\in \Irr(M)$, 
we say that  $$(N, \nu)\leq (M, \mu)$$
\noindent if $N\leq M$ and $[\nu, \mu_N]\neq 0$. Similarly  we say that
$(N, \nu)<(M, \mu)$ if 
$N<M$ and $[\nu, \mu_N]\neq 0$.

Let $N$ be a normal subgroup of $G$.
Let  $\chi\in \Irr(G)$ and $\nu\in \Irr(N)$ be
 such that
$(N, \nu)\leq (G, \chi)$. We denote by $G_{\nu}=\{g\in G\mid \nu^g=\nu\}$
the stabilizer of $\nu$ in $G$. Also we denote by
$\chi_{\nu}\in\Irr(G_{\nu})$ the Clifford correspondent of $\chi$ with respect
to $\nu$, i.e. the unique character $\chi_{\nu}\in \Irr(G_{\nu})$ 
such that
$(\chi_{\nu})_N =\frac{\chi_{\nu}(1)}{\nu(1)}\nu$ and $(\chi_{\nu})^G=\chi$.

 We set  $\Irr(M \modu N)=\{\mu \in \Irr(M)\mid \Ker(\mu)\geq N\}$
and $\Lin(G)=\{\chi\in \Irr(G) \mid \chi(1)=1\}$. The principal 
character of a group $H$ is denoted by $1_H$.  
\end{section}

\begin{section}{Preliminaries}
In this section, we will prove a series of lemmas that will be used 
in the proof of Theorems A and B.

\begin{hyp}\label{hypmain} Let $p$ be a prime number and $G$ be a finite 
$p$-group. Let $\chi \in \Irr(G)$ be such that $\chi(1)=p^n$. 
\end{hyp}
\begin{lem}\label{basic}
Assume Hypotheses \ref{hypmain}.
Let $H$ be a subgroup of $G$ and $\gamma\in \Irr(H)$.
 Assume that $\gamma^G=\chi$. Then 
 $\Xi( (\gamma\overline{\gamma})^G)\subseteq  \Xi (\chi\overline{\chi})$.
\end{lem}
\begin{proof} Since $\gamma^G=\chi$, we have that
$[\overline{\gamma}, \overline{\chi}_H]\neq 0$. 
By Exercise 5.3 of  \cite{isaacs} we have that 
$(\gamma \overline{\chi}_H)^G=\chi\overline{\chi}$. 
Thus $\chi\overline{\chi}= (\gamma\overline{\gamma})^G +\Theta$, for some
character $\Theta$ of $G$.
\end{proof}

\begin{lem}\label{basic2}
Assume Hypotheses \ref{hypmain}.
Let $M/N$ be a chief factor of $G$. Assume
 $\Xi ( (1_N)^M) \subseteq \Xi((\chi\overline{\chi})_M)$. Then for
each $\delta\in \Irr(M\modu N)$, there exists $\theta_{\delta} \in 
\Xi(\chi\overline{\chi})$ such that
 $(\theta_{\delta})_M=\theta_{\delta}(1) \delta$. In particular, 
there are at least $p-1$ distinct elements in 
$\Xi (\chi\overline{\chi})$ lying above
 irreducible constituents of $ ( (1_N)^M-1_M)$.  
\end{lem}
\begin{proof}
Since $G$ is a $p$-group, $G$ acts trivially on 
$M/N$. Therefore every $\delta \in \Irr(M\modu N)$ is 
$G$-invariant and linear. By hypothesis we have that, 
given $\delta \in \Irr(M\modu N)$, there exists some 
$\theta_{\delta} \in \Xi(\chi\overline{\chi})$ such that
 $[(\theta_{\delta})_M, \delta]\neq 0$. 
Since $\delta$ is $G$-invariant and linear, we have that
 $(\theta_{\delta})_M=\theta_{\delta} (1) \delta$. 
Thus if $\delta_1, \delta_2 \in \Irr(M\modu N)$ and 
$\delta_1\neq \delta_2$, then $\theta_{\delta_1}\neq \theta_{\delta_2}$.
\end{proof}

\begin{lem}\label{induction}
Assume Hypotheses \ref{hypmain}. 
Let  $M/N$ be a chief factor of $G$. Let $\mu\in \Irr(M)$ and
 $\nu\in \Irr(N)$
be characters such that  
$(N, \nu)< (M,\mu)\leq (G,\chi)$.
  Assume that $G_{\mu}< G_{\nu}$. Then  $|G_{\nu}: G_{\mu}|=p$ and
$\Xi((1_N)^M)\subseteq \Xi((\chi\overline{\chi})_{M})$.
\end{lem}
\begin{proof}
Let
$\chi_{\nu}\in \Irr(G_{\nu})$ be the Clifford
correspondent of $\chi$ with respect to $\nu$. Therefore $(\chi_{\nu})_N$
is a multiple of $\nu$. Because $G_{\mu}< G_{\nu}$ and $|M/N|=p$,
 there are $p$ conjugates
of $\mu$, and they are all the characters of $M$ lying above $\nu$. 
Therefore $|G_{\nu}: G_{\mu}|=p$
and $(\chi_{\nu})_M$ is a multiple of the sum of these
 $p$-conjugates. The
last sum is  $\nu^M$.
Therefore  $(\chi_{\mu})_M$ is a multiple of $\nu^M$.
Thus $\Xi(\nu^M)\subseteq \Xi((\chi_{\nu})_M)$. We conclude that
$\Xi((1_N)^M) \subseteq \Xi( (\nu\overline{\nu})^M)\subseteq
\Xi( \nu^M(\overline{\nu})^M ) 
\subseteq \Xi((\chi_{\nu}\overline{\chi_{\nu}})_M)\subseteq 
\Xi((\chi\overline{\chi})_{M})$.
\end{proof}

\begin{lem}\label{induccioncase} Assume Hypotheses \ref{hypmain}.
 Let $N\triangleleft G$
and $\nu\in \Irr(N)$ be such that $(N,\nu)< (G,\chi)$.
 Let $M/N$ be a chief factor of 
$G$. Assume that $M\cap G_{\nu}=N$. Set $\mu =\nu^M$. 
 Then $\mu\in \Irr(M)$,  $G_{\mu}=G_{\nu}M$ and $(M, \mu)\leq (G, \chi)$. 
Also $\Xi((1_N)^M)\subseteq \Xi( \chi\overline{\chi})$.
\end{lem}
\begin{proof}
Since  $G_{\nu}\cap M=N$,
 by Clifford theory it follows that $\mu=\nu^M\in \Irr(M)$, and that
$\mu_N$ is the sum of the $p$ distinct $M$-conjugates of $\nu$.

 Let $\chi_{\nu}\in \Irr(G_{\nu})$
be the Clifford correspondent of $\chi$ with respect to $\nu$.
By Mackey's Theorem we have that $((\chi_{\nu})^{(G_{\nu}M)})_M=
((\chi_{\nu})_N)^M$. Since $((\chi_{\nu})_N)^M$ is a multiple of 
$\mu$, it follows that  
$G_{\nu}M\leq G_{\mu}$. 
Since $\mu=\nu^M\in \Irr(M)$ is the unique character of $M$
above $\nu$ and $[\chi_N, \nu]\neq 0$, we have
that $[\chi_M, \mu]\neq 0$. Let $\chi_{\mu}\in \Irr(G_{\mu})$
be the Clifford correspondent of $\chi$ with respect to 
$\mu$. Observe that 
$$(\chi_{\mu})_N=((\chi_{\mu})_M)_N=( \frac{\chi(1)}{\mu(1)}  \mu)_N.$$
\noindent This last expression is a multiple of the $p$ distinct 
$M$-conjugates of $\nu$, which must therefore be the 
$G_{\mu}$-conjugates of $\nu$. Therefore $|G_{\mu}: G_{\nu}|=p$.
 Since $|G_{\nu}M: G_{\nu}|=p$
and $G_{\nu}M\leq G_{\mu}$, it follows
that $G_{\mu}=G_{\nu}M$. 

Since $\nu^M=\mu$, by Lemma \ref{basic} we have that
$\Xi( (\nu\overline{\nu})^M)\subseteq \Xi (\mu\overline{\mu})$.
Observe that $\Xi ((1_N)^M)\subseteq \Xi ((\nu\overline{\nu})^M)$
and $\Xi (\mu\overline{\mu})\subseteq  \Xi( (\chi\overline{\chi})_M)$.
Thus $\Xi((1_N)^M)\subseteq \Xi( \chi\overline{\chi})$.
\end{proof}

\begin{lem}\label{induc}
Assume Hypotheses \ref{hypmain}. Assume also that $\chi(1)\geq p$. 
Set $Z={\bf Z}(\chi)$. Then $Z<G$. 
 Let $M/Z$ be a chief factor of $G$. Let $\nu\in \Lin(Z)$
 be the unique 
character such that $\chi_{Z}=\chi(1) \nu$. Fix $\mu\in \Irr(M)$
such that $(Z, \nu)< (M, \mu)\leq (G, \chi)$.  
Then 

(i)  $\Xi(( 1_{Z})^{M})\leq \Xi((\chi\overline{\chi})_{M})$.

(ii) $\mu(1)= p^0=1$ and $|G:G_{\mu}|=p$.

(iii) $\Xi(( 1_{G_{\mu}})^{G})\leq \Xi(\chi\overline{\chi})$.
\end{lem}
\begin{proof}
 Without loss of 
generality, we may assume that $\Ker(\chi)=1$ and thus
$Z={\bf Z}(\chi)={\bf Z}(G)$. 
 Since 
$\chi(1)>1$, we have that $Z< G$. 
 Since $Z={\bf Z}(\chi)$
and $|M:Z|=p$, it follows that $(\mu)_Z=\nu$. Observe that 
$\chi_{M}\neq \chi(1)\mu$ since $M> Z={\bf Z}(\chi)$.
 Thus $G_{\mu}<G= G_{\nu}$.  Observe that (i) and (ii) follow 
by Lemma \ref{induction}.
Let  $\chi_{\mu}\in \Irr(G_{\mu})$ be the Clifford correspondent
of $\chi$ with respect to $\mu$. Then $(\chi_{\mu})^G=\chi$ and 
by Lemma \ref{basic} we have that 
$\Xi((\chi_{\mu}\overline{\chi_{\mu}})^G) \subseteq 
\Xi(\chi\overline{\chi})$. 
Therefore  $\Xi(( 1_{G_{\mu}})^{G})\subseteq
 \Xi((\chi_{\mu}\overline{\chi_{\mu}})^G)\subseteq 
\Xi(\chi\overline{\chi})$
and (iii) follows.
\end{proof}

\begin{lem}\label{induction2}
Assume Hypotheses \ref{hypmain}. 
Let 
$M/N$ be a chief factor of $G$.
Let $\nu\in \Irr(N)$ be such that $[\nu, \chi_N]\neq 0$.
Assume that $\nu^{M}=\mu\in \Irr(M)$. Then 
$\Xi((1_N)^M)\subseteq \Xi((\chi\overline{\chi})_{M})$.
\end{lem}
\begin{proof}
Observe that $[\mu, \chi_M]\neq 0$. By Lemma
\ref{basic} we have that
$\Xi( (\nu\overline{\nu})^M)\subseteq \Xi(\mu\overline{\mu})$.
Thus $
 \Xi((1_N)^M)\subseteq \Xi( (\nu\overline{\nu})^M)\subseteq 
\Xi(\mu\overline{\mu})$. 
Since $\Xi(\mu\overline{\mu})\subseteq  \Xi((\chi\overline{\chi})_{M})$,
the result follows.
\end{proof}
\end{section}
\begin{section}{Theorem A}
\begin{hyp}\label{hyp2}
 Assume Hypotheses \ref{hypmain}. Let
$1=N_0< N_1<\cdots<N_j<\cdots < N_t=G$
be a composition series of $G$. Pick $\nu_j\in \Irr(N_j)$, 
for $j=0,1,\ldots, t$, such that 
\begin{equation}\label{chain}
(1, 1_1)=
(N_0, \nu_0)<(N_1, \nu_1)<
\cdots< ( N_j, \nu_j)< \cdots<( N_t, \nu_t)=
(G,  \chi). 
\end{equation}
\end{hyp}
\begin{definition} Assume Hypotheses \ref{hyp2}. Fix $i=1, \ldots, t$.
We say that the chain \eqref{chain}
is {\bf stable at $i$} if $\nu_{i-1}=(\nu_i)_{N_{i-1}}$
and $G_{\nu_{i-1}}=G_{\nu_i}$. 
Otherwise we say that that chain is {\bf unstable at $i$}.
\end{definition}

\begin{lem}\label{unstable} Assume Hypotheses \ref{hyp2}. Fix
$i=1, \ldots,  t$. 
If the chain \eqref{chain} is unstable at $i$, then
 $\Xi((1_{N_{i-1}})^{N_i})\subseteq \Xi( \chi\overline{\chi})$.
Also
one of the following holds:
	
(a) $\nu_{i-1}=(\nu_i)_{N_{i-1}}$,  $G_{\nu_{i}}<G_{\nu_{i-1}}$
and $|G_{\nu_{i-1}}: G_{\nu_{i}}|=p$.

(b) $\nu_i=(\nu_{i-1})^{N_i}$ and
$G_{\nu_i}=G_{\nu_{i-1}}N_i$. Thus $|G_{\nu_i}:G_{\nu_{i-1}}|=p$.
\end{lem}

\begin{proof} Since $N_{i}/N_{i-1}$ is a chief
factor of a $p$-group, it is cyclic of order $p$.  Either 
$\nu_{i-1}$ extends to $N_i$, or $(\nu_{i-1})^{N_{i}}\in \Irr(N_i)$.

Assume that $\nu_{i-1}$ extends to $N_i$. Since $N_{i}/N_{i-1}$
is a cyclic group, it follows that
$\nu_{i-1}=(\nu_i)_{N_{i-1}}$. Observe that  
$(\nu_i)^g=\nu_i$ implies that $(\nu_{i-1})^g=\nu_{i-1}$ for any $g$ in $G$.
Thus  $G_{\nu_{i}}\leq G_{\nu_{i-1}}$. Since the chain is 
unstable at $i$,  we have that  $G_{\nu_{i}}<G_{\nu_{i-1}}$.
By Lemma \ref{induction}
 we have that 
  $|G_{\nu_{i-1}}: G_{\nu_i}|=p$
and $\Xi((1_{N_{i-1}})^{N_i})\subseteq \Xi((\chi\overline{\chi})_{N_i})$.
Hence the lemma holds with option (a).
  
We may assume now that $(\nu_{i-1})^{N_{i}}\in \Irr(N_i)$.
Since $[(\nu_i)_{N_{i-1}}, \nu_{i-1}]\neq 0$, we have that 
that $\nu_i= (\nu_{i-1})^{N_{i}}$. 
Also $G_{\nu_{i-1}}\cap N_i=N_{i-1}$. By Lemma 
\ref{induccioncase} we have that 
 $G_{\nu_i}=G_{\nu_{i-1}}N_i$  and
 $\Xi((1_{N_{i-1}})^{N_i})\subseteq \Xi( \chi\overline{\chi})$.
So $|G_{\nu_i}:G_{\nu_{i-1}}|=|N_i: G_{\nu_{i-1}}\cap N_i|=|N_i:N_{i-1}|=p$.
Hence the lemma holds with option (b) and the proof is complete.
\end{proof}
\begin{lem}\label{basicta}
Assume Hypotheses \ref{hyp2}.      
For $i=1, \ldots ,t$,
let $r_i$ and $s_i$ be such that $p^{r_i}=|G:G_{\nu_i}|$ and $
p^{s_i}= \nu_i(1)$. Set
\begin{equation}\label{mi}
m_i=
|\{k\mid 1\leq k\leq i  
\mbox{ and the chain \eqref{chain} is unstable at $k$ }\}|.
\end{equation}

Then for any $i=1,\ldots, t$

{\bf (i)}  $\Xi ( (1_{N_{i-1}})^{N_i}) \subseteq
 \Xi((\chi\overline{\chi})_{N_i})$
if the chain is unstable at $i$.

{\bf (ii)}  $m_i= 2s_i + r_i$.
\end{lem}
\begin{proof}
By Lemma \ref{unstable} we have that (i) holds. It remains to 
prove that (ii) also holds.

Since $G$ is a $p$-group and $N_1/1$ is a chief factor of $G$, 
we have that  $N_1\leq {\bf Z}(G)$. Therefore 
 $p^{r_1}=|G:G_{\nu_1}|=|G:G|=1$  and $p^{s_1}=\nu_1(1)=1$.
Thus $r_1=0$ and $s_1=0$. Since
$i$ is stable, we have that $m_1=0$. So (ii) holds for $i=1$. 

We may assume by induction that we have some integer $j=2,\ldots, t$
such that the lemma holds for any 
$i\leq j-1$. 
Either the chain is stable at $j$, or it is unstable
at $j$.

If the chain is stable at $j$, then by \eqref{mi}
we have that $m_j=m_{j-1}$.
Also  $p^{r_j}=|G:G_{\nu_j}|=|G:G_{\nu_{j-1}}|
=p^{r_{j-1}}$  and $p^{s_j}=\nu_j(1)=\nu_{j-1}(1)= p^{s_{j-1}}$.
Thus $$m_j = m_{j-1}= 2s_{j-1} + r_{j-1}= 2s_j + r_j.$$
\noindent So (ii) holds for $i=j$.

Now we may assume that the chain is unstable at $j$. By
\eqref{mi} we have that $m_j = m_{j-1}+1$. Also
 either (a) or (b)  of Lemma \ref{unstable} holds.
If  (a) holds, then 
  $p^{r_j}=|G:G_{\nu_j}|=p|G:G_{\nu_{j-1}}|
=p^{r_{j-1}+1}$  and $p^{s_j}=\nu_j(1)=\nu_{j-1}(1)= p^{s_{j-1}}$. 
Thus $$m_j = m_{j-1}+1 = 2s_{j-1} + r_{j-1}+ 1= 2s_j + r_j.$$
\noindent So (ii) holds.

Finally,  assume that Lemma \ref{unstable} (b) holds. 
Thus
$p p^{r_j}=p |G:G_{\nu_j}|=|G:G_{\nu_{j-1}}|=p^{r_{j-1}}$. 
Therefore
$r_j= r_{j-1}-1$.  Since  $p^{s_j}=\nu_j(1)=p\nu_{j-1}(1)= p p^{s_{j-1}}$, 
we have that $s_j=s_{j-1}+1$ and 
$$m_j = m_{j-1}+1 = 2s_{j-1} + r_{j-1}+ 1= (2s_{j-1} +2) + (r_{j-1}-1) =
2s_j +r_j.
$$
\noindent So (ii) holds and the proof is complete.
\end{proof}

\begin{proof}[Proof of Theorem A]
We may assume that Hypotheses \ref{hyp2} hold.

Fix $i=1,\ldots, t$.
Let $\theta\in \Xi(\chi\overline{\chi})$ be a character
that lies 
above  an irreducible constituent of $((1_{N_{i-1}})^{N_i}-1_{N_i})$.
If $j>i$, then $[\theta_{N_j},\delta]=0$ for any 
$\delta \in \Irr(N_j\modu N_{j-1})$. Otherwise 
$\Ker(\delta)\geq N_{j-1}\geq N_i$. Therefore
$[1_{N_{j-1}}, \theta_{N_{j-1}}]\neq 0$. Since
$N_{j-1}\trianglelefteq G$, we have that $N_{j-1}\leq \Ker(\theta)$.
Thus  $\theta_{N_{i}}=\theta(1)1_{N_i}$ and $\theta$ does not
lie above an irreducible constituent of $((1_{N_{i-1}})^{N_i}-1_{N_i})$.
If $j<i$ then $\theta_{N_j}=\theta (1) 1_{N_j}$.
We conclude that for any $i, j=1, \ldots t$ with $i\neq j$,  the elements of 
$\Xi(\chi\overline{\chi})$ lying
above irreducible constituents of $((1_{N_{i-1}})^{N_i}-1_{N_i})$
are distinct from the elements of 
$\Xi(\chi\overline{\chi})$ lying
above  irreducible constituents of $((1_{N_{j-1}})^{N_j}-1_{N_j})$.

Using the notation of Lemma \ref{basicta}, we have that
$s_t=n$ and $r_t=0$. Thus $m_t=2n$ by Lemma \ref{basicta} (ii).
 If the chain is unstable at $i=1,\ldots, t$, then 
by Lemma \ref{basic2} and Lemma \ref{basicta} (i) we have that
there are at least $p-1$ non-principal distinct 
 elements of $\Xi(\chi\overline{\chi})$ lying
above  irreducible constituents of $((1_{N_{i-1}})^{N_i}-1_{N_i})$.
Therefore, by the previous paragraph
there are at least $m_t(p-1)$ distinct non-principal
elements of $\Xi(\chi\overline{\chi})$. Since $1_G$ is an element of 
$\Xi(\chi\overline{\chi})$, we conclude that 
$\Xi(\chi\overline{\chi})$ has at least $m_t(p-1)+1= 2n(p-1)+1$
distinct elements.
\end{proof}

A corollary of Theorem A is 

\begin{cora}\label{corola}
Let $p$ be a prime number. Let $G$ be a finite $p$-group and 
$\chi, \psi \in \Irr(G)$ with $\chi(1)=p^n$. Assume that
there exists $\alpha \in \Lin(G)$ such that
$[\alpha, \chi\psi]\neq 0$. Then 
$$\eta(\chi, \psi)\geq 2n(p-1)+1.$$
\end{cora}
\begin{proof}
Observe that
$[\chi\psi, \alpha]= [\psi, \alpha \overline{\chi} ]\neq 0$. 
Since $\alpha(1)=1$, we have that
$\psi=\alpha\overline{\chi}$. 

Let $\theta_1, \theta_2 \in \Irr(G)$.
Observe that $\alpha \theta_1=\alpha\theta_2$ if and only 
if $\theta_1=\theta_2$. Also observe that $[\theta_1,\chi\overline{\chi}]=
[\alpha\theta_1, \chi\psi]$. Since $\alpha\theta\in \Irr(G)$
for any $\theta\in \Irr(G)$, it follows that 
$\eta(\chi, \psi)=\eta(\chi, \overline{\chi})$. Applying 
Theorem A, we get that $\eta(\chi, \psi)\geq 2n(p-1)+1$.
\end{proof}

 Let $p$ be an odd prime and $n>0$ be an integer.  
Let $E$ be an extra-special group of order $p^{2n+1}$ and exponent 
$p$. Let $\chi\in \Irr(E)$ be
such that $\chi(1)=p^n$. Since $p>2$, we can choose 
$\psi\in \Irr(E)$ such that
$\psi(1)=p^n$ and $\psi\neq \overline{\chi}$. We can check that 
$\eta(\chi, \psi)=1$. Thus Corollary \ref{corola} does not
remain true without the existence  of $\alpha\in \Lin(G)$ such 
that $[\chi\psi, \alpha]\neq 0$. 

\end{section}

\begin{section}{Theorem B}

\begin{lem}\label{casep}
 Assume Hypotheses \ref{hypmain}.
Assume also that $\chi(1)=p$.
Let  $\theta\in \Irr(G)$
such that $[\theta, \chi\overline{\chi}]\neq 0$.
 Then $[\theta, \chi\overline{\chi}]=1$
and
one of the following holds:

(i) $\eta(\chi, \overline{\chi})=2p-1$ and $\chi\overline{\chi}$ is
the sum of $p$ distinct linear characters and $p-1$ distinct 
irreducible characters
of degree $p$.

(ii) $\eta(\chi, \overline{\chi})=p^2$ and $\chi\overline{\chi}$
is the sum of $p^2$ distinct linear characters.
\end{lem}
\begin{proof}
Without loss of generality, we may assume that $\Ker(\chi)=1$. 
Set $Z={\bf Z}(G)$.
Let $M/Z$ be a chief factor of $G$. Let $\nu\in\Lin(Z)$ be 
the unique character
of $Z$ such that $\chi_Z=\chi(1)\nu$.  Let $\mu\in \Irr(M)$ be
such that $(Z,\nu)<(M,\mu)\leq (G, \chi)$.
Set $H=G_{\mu}$.
Let $\chi_{\mu}\in \Irr(H)$ be the Clifford correspondent of $\chi$
with respect to $\mu$. By Lemma \ref{induc} (ii) we have that $|G:H|=p$.

 Since $\chi(1)=p$, $|G:H|=p$, $\chi_{\mu}\in \Irr(H)$
and $(\chi_{\mu})^G=\chi$,
we have that $\chi_{\mu}(1)=1$.
By Exercise 2.8 of \cite{isaacs}
we have that $H$ is an abelian subgroup of $G$
since $\chi$ is faithful and $H$ is normal. 

	Fix $g\in G \setminus H$. Since $H$ is normal and 
$(\chi_{\mu})^G=\chi$, we have that 
\begin{equation}\label{equationchi}
	(\chi \overline{\chi})_H= \sum_{i=0}^{p-1} (\chi_{\mu})^{g^i}
 \sum_{j=0}^{p-1} (\overline{\chi_{\mu}})^{g^j}=
 \sum_{i=0}^{p-1}\sum_{j=0}^{p-1} [\chi_{\mu}
 (\overline{\chi_{\mu}})^{g^j}]^{g^i}.
\end{equation}

	Consider the set 
$\Omega=\{\chi_{\mu}(\overline{\chi_{\mu}})^{g^j}\mid j=0,1, \ldots , p-1\}$. 
Observe that
the elements of $\Omega$ are linear characters of $H$ and $\Omega$ has
exactly $p$ elements.
Also observe that, given $i, j=0, 1,\ldots, p-1$, the product
$(\chi_{\mu})^{g^i} (\overline{\chi_{\mu}})^{g^j} \in \Irr(H)$
 is $G$-conjugate to an element in $\Omega$. 

\begin{claim}\label{cla}
 Assume that for some $s\in \{1, \ldots , p-1\}$, the character
$\chi_{\mu}(\overline{\chi_{\mu}})^{g^s}$ is $G$-invariant.
Then for all $i=0, 1, \ldots, p-1$, the character
$\chi_{\mu}(\overline{\chi_{\mu}})^{g^i}$ is $G$-invariant.
Thus the set $\Xi((\chi_{\mu}(\overline{\chi_{\mu}})^{g^i})^G)$ contains $p$
linear characters. Also
 $\Xi((\chi_{\mu}(\overline{\chi_{\mu}})^{g^i})^G)$ and 
$\Xi((\chi_{\mu}(\overline{\chi_{\mu}})^{g^j})^G)$ are disjoint
sets for all $i,j=0,1,\ldots, p-1$ with $i\neq j$. In particular,
the set $\Xi(\chi\overline{\chi})$ contains  $p^2$ 
distinct elements, all of which are linear characters.
Thus $\eta(\chi, \overline{\chi})=p^2$ and 
Lemma \ref{casep} (ii) holds.
\end{claim}
  
\begin{proof}
Since $\chi_{\mu}(\overline{\chi_{\mu}})^{g^s}$ is $G$-invariant,
 $K=\Ker(\chi_{\mu}(\overline{\chi_{\mu}})^{g^s})$ is a normal
subgroup of $G$ with $Z\leq K\leq H$. 
	
Suppose that $K>Z$. We choose $L\leq K$ such that
$L/Z$ is a chief factor of $G$. Let $\beta= (\chi_{\mu})_L$.
Observe that $\beta\in \Lin(L)$ since $\chi_{\mu}(1)=1$.
If $G_{\beta}=G$, then $\chi_L=\chi(1)\beta$ and therefore
$L\leq {\bf Z}(\chi)=Z$. Since $L/Z$ is a chief factor of $G$,
we conclude that
 $G_{\beta}=H$.  Observe that 
$(\chi_{\mu})^{g^s}_L= \beta^{g^s}\neq \beta$ since $g^s \in G\setminus H$. 
Thus $L \not \leq K$ since 
$\beta \overline{\beta^{g^s}}\neq 1_L$. 
We conclude that $K=Z$.

Now  $\chi_{\mu}(\overline{\chi_{\mu}})^{g^s}$ is a $G$-invariant linear
character with kernel $K=Z$. So $G$ centralizes $H/Z$. Hence 
all the characters in $\Irr (H \modu Z)$ are $G$-invariant.
 Since $\Ker(\chi\overline{\chi})=Z$, we have that
$\Ker( \chi_{\mu}(\overline{\chi_{\mu}})^{g^i})\geq Z$
            for all $i=0, 1, \ldots, p-1$.
Thus  the character
$\chi_{\mu}(\overline{\chi_{\mu}})^{g^i}\in \Lin(H)$ is $G$-invariant 
 for all $i=0, 1, \ldots, p-1$.

 Since   $\chi_{\mu}(\overline{\chi_{\mu}})^{g^i}$ is $G$-invariant,
for $i=0,1, \ldots, p-1$,
and $|G:H|=p$, it follows that  $\chi_{\mu}(\overline{\chi_{\mu}})^{g^i}$
extends to $G$ and $(\chi_{\mu}(\overline{\chi_{\mu}})^{g^i})^G$ is the sum
of the $p$ distinct extensions of $\chi_{\mu}(\overline{\chi_{\mu}})^{g^i}$.
Since $\chi_{\mu}(\overline{\chi_{\mu}})^{g^i}\neq
 \chi_{\mu}(\overline{\chi_{\mu}})^{g^j}$ for any $i,j=0,1,\ldots, p-1$
with $i\neq j$, the set of extensions of $\chi_{\mu}(\overline{\chi_{\mu}})^{g^i}$
is disjoint from 
the set of extensions of $\chi_{\mu}(\overline{\chi_{\mu}})^{g^j}$.
Since $((\chi_{\mu})^{g^i})^G=\chi$ for any $i=0,1,\ldots,p-1$, we have that
 $\Xi((\chi_{\mu}(\overline{\chi_{\mu}})^{g^i})^G)\subseteq  
\Xi((\chi_{\mu})^G((\overline{\chi_{\mu}})^{g^i})^G)=\Xi(\chi\overline{\chi})$.
We conclude that the set $\Xi(\chi\overline{\chi})$ has at least
$p\times p=p^2$ distinct elements. Since $\chi(1)=p$, it follows
that $\eta(\chi, \overline{\chi})=p^2$.
\end{proof}

If for some $s=1,\ldots,p-1$, the character
$\chi_{\mu}(\overline{\chi_{\mu}})^{g^s}$ is $G$-invariant, the lemma
follows by Claim \ref{cla}.

 Now we assume that none of the characters
$\chi_{\mu}(\overline{\chi_{\mu}})^{g^i}$, for $i=1,\ldots, p-1$, 
is $G$-invariant. Since $|G:H|=p$, it follows that
$(\chi_{\mu}(\overline{\chi_{\mu}})^{g^i})^G \in \Irr(G)$
and $(\chi_{\mu}(\overline{\chi_{\mu}})^{g^i})^G(1)=p$ for
$i=1,\ldots, p-1$.
By \eqref{equationchi} we have $[(\chi\overline{\chi})_H,
\chi_{\mu}(\overline{\chi_{\mu}})^{g^i})]\neq 0$. Thus
$$[\chi\overline{\chi},
(\chi_{\mu}(\overline{\chi_{\mu}})^{g^i})^G]\neq 0.$$
We conclude that the set 
$\{((\chi_{\mu})^{g^j}(\overline{\chi_{\mu}})^{g^i})^G
\mid i,j=0,1,\ldots,p-1 \mbox{ and } i\neq j\}$ has at most $p-1$ 
characters, all of which are irreducible
characters of degree $p$, and lying in the set 
$\Xi(\chi\overline{\chi})$.

Since ${\chi_{\mu}}^G=\chi$, by Lemma \ref{induccioncase} we have that
$\Xi((1_H)^G)\subseteq \Xi(\chi\overline{\chi})$.
Observe that  $(1_H)^G$ has $p$ distinct irreducible constituents.
Observe also that $\Xi((1_H)^G)$ and 
$\{((\chi_{\mu})^{g^j}(\overline{\chi_{\mu}})^{g^i})^G
\mid i,j=0,1,\ldots,p-1 \mbox{ and } i\neq j\}$ are disjoint
sets since $(\chi_{\mu})^{g^j}(\overline{\chi_{\mu}})^{g^i}\neq 1_H$
for any $i,j=0,1,\ldots,p-1$ with $i\neq j$.
We conclude that the set $\Xi(\chi \overline{\chi})$ has at
most $2p-1$ members, i.e. $\eta(\chi,\overline{\chi})\leq 2p-1$. By 
Theorem A we have that 
$\eta(\chi, \overline{\chi})\geq 2p-1$ and therefore
 $\eta(\chi, \overline{\chi})= 2p-1$.
 Since $(p-1)p+ p=p^2$,
we have that all the irreducible constituents of $\chi\overline{\chi}$
appear with multiplicity $1$ and have the right degree.
So the proof is complete.
\end{proof}

\begin{proof}[Proof of Theorem B]
Theorem B follows from Lemma \ref{casep} and Theorem A. 
\end{proof}
\end{section}

\begin{section}{Examples}
\begin{prop}\label{example}
Let $p$ be a prime number and $n\geq 0$ be an integer. 
There exist a finite
$p$-group $G$ and a character $\chi \in \Irr(G)$ such that 
$\chi\neq \overline{\chi}$,
$\eta(\chi, \overline{\chi})=2n(p-1)+1$ and $\chi(1)= p^n$.
\end{prop} 
\begin{proof}
The statement is clearly true for $n=0$. So we may assume
that $n>0$.
By induction on $n$, we may assume that there exist a finite $p$-group
$A$ and a character $\alpha\in \Irr(A)$ such that 
$\alpha\neq \overline{\alpha}$, 
$\eta(\alpha, \overline{\alpha})=2(n-1)(p-1)+1$ and $\alpha(1)=p^{n-1}$.

Let $C=\{0,1, \ldots, p-1\}$ be the additive group of
 integers
modulo $p$. Let
 $H=A\times\cdots\times A$  be the direct product of $p$ copies of 
$A$. Let $\gamma$ be a generator of $C$. Define
the action of $\gamma$ in $H$ by
$$(a_0,a_1,\ldots,a_{p-1})^{\gamma}= (a_{p-1}, a_0, \ldots , a_{p-2})$$
\noindent for any $a_0, a_1, \ldots, a_{p-1}\in A$.

Let $G$ be the semi-direct product of $H$ and $C$.
Then  $G$ is  the wreath product $A\wr C$ of $A$ by $C$. 

For $i=0, \dots, p-1$, set $\theta_i((a_0,a_1, \ldots, a_{p-1}))= \alpha(a_i)$.
Observe that $\theta_i\in \Irr(H)$.
Also observe that the stabilizer of $\theta_0$ is $H$.
Thus $\chi =\theta_0^G \in \Irr(G)$ and 
$\chi(1)=p \theta_0(1)=p\alpha(1)=p^n$.

\begin{claim} $\chi\neq \overline{\chi}$.
\end{claim}
\begin{proof} 	
Since $\alpha\neq \overline{\alpha}$, there exists some $a\in A$
such that $\alpha(a)\neq \overline{\alpha(a)}$. 
Observe that 
\begin{equation*}
\begin{split}
\chi((a,1,\ldots,1))&= \theta_0((a,1,\ldots,1))+ \theta_0((1, a,\ldots,1)+\cdots 
+\theta_0((1,1,\ldots,a))\\
&= \alpha(a)+\alpha(1)+ \cdots +\alpha(1)\\
& =\alpha(a) +(p-1)\alpha(1)\\
&\neq
\overline{\alpha(a) +(p-1)\alpha(1)}\\
&=\overline{\chi((a,1,\ldots,1))},
\end{split}
\end{equation*}
\noindent and the claim follows.
\end{proof}
Observe that
$\chi_H= \sum_{i=0}^{p-1} \theta_i$. Consider the set 
 $\Omega= \{\theta_0 \overline{\theta_k}\mid k=1, \ldots, p-1\}$.
Observe that for $i, j = 0,1,\ldots ,p-1$, $i\neq j$,
we have that $\theta_i \overline{\theta_j} \in \Irr(H)$ and 
 $\theta_i \overline{\theta_j}$ is $G$-conjugate to 
the character $\theta_0 \overline{\theta_k}$, where $k\cong j-i \modu p$.
Since $\theta_0 \overline{\theta_i}\neq \theta_0 \overline{\theta_j}$
if $i\neq j$, the set $\Omega$ has $p-1$ distinct elements.

 For any $i =1,\ldots, p-1$, the character
 $\theta_0\overline{\theta_i}$ is not $G$-invariant.
Otherwise $\theta_0\overline{\theta_i}
=\alpha\times 1_A \times\cdots\times 1_A \times \overline{\alpha}\times 1_A\times \cdots \times 1_A$,
where $\overline{\alpha}$ is in the $i$ position, is equal
to 
$\overline{\theta_0} \theta_{p-i}
=\overline{\alpha}\times 1_A \times \cdots \times 1_A
\times\alpha\times 1_A \times \cdots\times
1_A$, where $\alpha$ is in the $p-i$ position. Hence
$\alpha=\overline{\alpha}$, a contradiction with our hypothesis.
Since $|G:H|=p$ and $\theta_0\overline{\theta_i}\in \Lin(H)$ is not 
$G$-invariant,
for $i=1,\ldots, p-1$, it follows that
 $(\theta_0 \overline{\theta_i})^G\in \Irr(G)$.
Therefore the set  $\{(\theta_0 \overline{\theta_i})^G\mid
i=1, \ldots, p-1\}= \{(\theta_i \overline{\theta_j})^G\mid
i,j =0, 1, \ldots, p-1 \mbox{ and }i\neq j \}$ has exactly $p-1$ distinct
irreducible  characters. 

Since 
$\theta_0\overline{\theta_0}=\alpha\overline{\alpha}
\times 1_A\times\cdots\times 1_A$
and  $\eta(\alpha, \overline{\alpha})= 2(n-1)(p-1)+1$,
it follows that  $\theta_0\overline{\theta_0}$ has $2(n-1)(p-1)+1$
distinct irreducible constituents. Among those constituents,
only the character $1_H$ is $G$-invariant. 
It follows that $(\theta_0\overline{\theta_0})^G$ has exactly 
 $p + 2(n-1)(p-1)$ distinct
irreducible constituents.  

We conclude that $\eta(\chi, \overline{\chi})=p + 2(n-1)(p-1) + p-1=
2n(p-1)+1$. 
\end{proof}

\end{section}    
{\bf Acknowledgment.} I would like to thank Professor Everett C. Dade 
for helpful discussions.

\end{document}